\documentclass[a4paper,12pt]{article}

\usepackage[french, american]{babel} 
\usepackage[utf8]{inputenc}
\usepackage[a4paper, left=2.5cm, right=2.5cm, top=2.5cm, bottom=2.5cm]{geometry}
\usepackage{amssymb,amsmath,amsthm}
\usepackage{mathrsfs}
\usepackage{dsfont}
\usepackage{float}

\usepackage{graphicx}
\usepackage{caption}
\usepackage{color}
\usepackage{tikz}

\numberwithin{equation}{section} 
\usepackage{comment}
\usepackage{hyperref}

\newcommand{\la}{\lambda}

\renewcommand{\phi}{\varphi}

\newcommand{\dd}{\mathrm{d}}

\newcommand{\Sp}{\mathrm{Sp}}
\newcommand{\tr}{\mathrm{tr}}

\newcommand{\ri}{\mathrm{i}}

\renewcommand{\Im}{\mathrm{Im}}
\renewcommand{\Re}{\mathrm{Re}}

\newcommand{\W}{\mathscr{W}}

\newcommand{\disteq}{\stackrel{d}{=}}
\newcommand{\distconv}{\xrightarrow[N \rightarrow \infty]{d}}

\newtheorem{theorem}{Theorem}[section]

\newtheoremstyle{myrem}
  {}
  {}
  {}
  {}
  {\bfseries}
  {.}
  { }
  {\thmname{#1}\thmnumber{ #2}\thmnote{\normalfont{ (#3)}}}
\theoremstyle{myrem}

\begin{document}
\title{\vspace{-2cm} Une brève histoire des perturbations non-hermitiennes de rang un\footnote{Version préliminaire («version auteur») de l'article du même nom publié dans la \href{https://smf.emath.fr/publications/la-gazette-de-la-societe-mathematique-de-france-182-octobre-2024}{Gazette SMF 182}.}}
\author{Guillaume Dubach\\
CMLS, École polytechnique, 91120 Palaiseau, France\\
 guillaume.dubach$@$polytechnique.edu
\and 
Jana Reker\\
UMPA, ÉNS de Lyon, 46 Allée d'Italie, 69007 Lyon, France\\
jana.reker$@$ens-lyon.fr}

\maketitle

\begin{abstract}
Les perturbations de faible rang de matrices aléatoires ont été au c\oe ur de nombreux travaux ces vingt dernières années. En particulier, les cas non-hermitiens, moins représentés dans la littérature en règle générale, font ici l'objet d'une attention spéciale en raison de leurs applications à la physique et à l'étude des réseaux de neurones. Petit tour d'horizon.
\end{abstract}

\section{Introduction}
C'est une idée vieille comme le monde -- bon, n'exagérons pas: vieille comme l'algèbre linéaire. On considère une matrice $M\in \mathcal{M}_N(\mathbb{C})$, 
que l'on va \textit{perturber} (par exemple, additivement) par une autre matrice $R\in \mathcal{M}_N(\mathbb{C})$. On s'intéresse alors au spectre de la matrice $M+R$, que l'on compare
à ceux de $M$ et $R$. Cette idée de base admet de multiples déclinaisons. Dans ce qui suit, néanmoins, on se limitera aux cas suivants: tout d'abord, la perturbation $R$ sera \textit{de rang un} -- ce qui, il faut bien l'avouer, simplifie considérablement l'analyse\footnote{Sous les dehors trompeurs d'une excessive simplicité, ces perturbations de rang un, comme on va le voir, ont déjà une structure très riche; elles ont fait récemment l'objet d'un savoureux article de survey de Peter Forrester~\cite{Forrester_Review}.}. Ensuite, bien qu'une très vaste littérature concerne des situations où $M$ et $R$ sont normales (par exemple toutes deux hermitiennes, comme dans le célèbre phénomène de transition BBP \footnote{Pour Baik, Ben Arous et Péché, en référence à l'article pionnier \cite{BBP}.}), on se penche ici sur des cas \textit{anormaux}, à savoir:
\begin{enumerate}
\item le cas d'une matrice aux entrées indépendantes et identiquement distribuées (i.i.d.) sans symétrie particulière; 
\item le cas d'une matrice hermitienne perturbée de façon anti-hermitienne (modèle $H+\ri \Gamma$);
\item enfin, le cas d'une matrice unitaire, perturbée multiplicativement par une matrice non unitaire (modèle $UA$).
\end{enumerate}
Un point commun entre ces différentes situations est que $M$ sera toujours, pour nous, une matrice aléatoire, et que nous ferons varier continûment la perturbation $R$ afin de considérer la dynamique des valeurs propres, qui suivent des trajectoires aléatoires dans le plan complexe. 

\bigskip
Le phénomène principal qui nous intéresse est l'émergence d'un ou plusieurs \textit{outliers}. Ce mot d'anglais n'ayant pas encore d'équivalent exact en français, nous l'emploierons tel quel; mais il mérite que l'on en précise d'abord la définition: est appelée \textit{outlier}, en principe, toute valeur propre éloignée du reste du spectre, ou éloignée d'une zone considérée comme typique pour un spectre aléatoire. On s'attend par ailleurs, dans la plupart des cas, à ce que cette valeur propre soit localisable (il s'agira, par exemple, d'en borner les fluctuations). Plus précisément, on dira qu'une valeur propre $\lambda_1$ est un \textit{outlier fortement séparé} s'il existe deux domaines $D_{1}$ et $D_{2}$ tels que, avec \textit{forte probabilité},
\begin{enumerate}
\item[(i)] $\lambda_{1}\in D_{1}$;
\item[(ii)] $\lambda_{2},\cdots,\lambda_{N}\in D_{2}$;
\item[(iii)] et tel que $D_{1}$ soit \textit{loin} de $D_{2}$.
\end{enumerate}
Les domaines $D_1, D_2$ peuvent être aléatoires ou déterministes. Ils dépendent également de la taille $N$ des matrices -- et l'on s'intéresse généralement au régime dans lequel cette taille tend vers l'infini. Dans ce contexte, la notion de forte probabilité concernera des événements dont le complémentaire est de probabilité plus petite que n'importe quelle puissance négative de $N$, et la séparation entre les domaines $D_1$ et $D_2$ jouera sur différentes échelles en fonction de $N$.

\bigskip
Dans chacun des modèles présentés ci-dessous, la localisation du ou des outliers est une question centrale en vue des applications -- et c'est une bonne nouvelle, parce que c'est aussi ce à quoi les techniques disponibles nous donnent accès le plus directement. Mais il y a d'autres questions naturelles qui sont plus difficiles à aborder. L'étude précise des fluctuations de l'outlier, de la forme des trajectoires, ou encore de l'origine de l'outlier dans le spectre de $M$, nécessite une analyse plus poussée. Et puisque dans le monde non-hermitien, il convient de distinguer \textit{deux} types de vecteurs propres, c'est à l'étude d'une dynamique couplée entre valeurs propres, vecteurs propres à droite, et vecteurs propres à gauche, que ces questions nous invitent.

\section{Le labyrinthe de Ginibre: perturbations de matrices aléatoires non-hermitiennes}

Évoquons tout d'abord le cas d'une matrice aléatoire $M$ sans symétrie particulière. L'intérêt pour les perturbations de rang un d'une telle matrice a été suscité par l'article de Rajan et Abbott de 2006 \cite{RajanAbbott2006}, qui établit la pertinence de ce modèle pour l'étude des réseaux de neurones\footnote{Il s'agit, en l'occurrence, de réseaux de \textit{vrais} neurones; mais ces applications se situent à l'interface entre les neurosciences computationnelles, qui visent à comprendre le fonctionnement du cerveau, et l'étude de réseaux de neurones artificiels en vue de l'IA. Une autre application naturelle des perturbations de faible rang à ces domaines est présentée dans l'article récent \cite{Schuessler}.}. Le point de départ de Rajan et Abbott est la \textit{loi de Dale}, qui stipule qu'un neurone donné est soit inhibiteur soit excitateur de tous les autres neurones. Cela se traduit dans leur modèle par une matrice de rang un, de la manière suivante: pour modéliser la \textit{matrice synaptique} (qui décrit les interactions entre neurones), on part d'une matrice $R$ de valeurs moyennes, où les entrées des colonnes correspondant aux neurones excitateurs sont toutes à une certaine valeur positive $\mu_E>0$, et celles des colonnes correspondant aux neurones inhibiteurs à une valeur négative $\mu_I<0$, ce qui forme bien une matrice de rang un, que l'on perturbe ensuite par un bruit aléatoire $G$ (aux entrées essentiellement indépendantes) qui correspond aux fluctuations attendues autour des valeurs moyennes. On souhaite comprendre la distribution des valeurs propres de $G+R$, modèle idéalisé de matrice synaptique. Cela rejoint bien sûr la philosophie générale exposée plus haut, à un détail près -- non pas mathématique mais linguistique: pour nous, c'est la matrice $R$ qui est considérée comme une `perturbation' de faible rang de la matrice aléatoire $G$. Formellement, cette différence de point de vue ne change rien; on cherche de toutes façons à considérer ce système dans tous les régimes possibles: on étudiera donc, par exemple, $G+tR$, avec $t$ un paramètre que l'on peut soit fixer soit faire varier continûment.

\bigskip
Les questions soulevées par Rajan et Abbott \cite{RajanAbbott2006} ont été évoquées par Terence Tao lors d'un workshop à l'American Institute of Mathematics, auquel participaient également Alice Guionnet et Percy Deift. Alice Guionnet avait récemment travaillé, avec Mylène Maïda et Florent Benaych-Georges, sur un problème analogue dans le cas hermitien \cite{BenaychGuionnetMaida}; Percy Deift quant à lui aurait appuyé l'idée d'utiliser l'identité de Sylvester, à savoir, que pour tout $A \in \mathcal{M}_{n \times d}(\mathbb{C})$ et $B \in \mathcal{M}_{d \times n} (\mathbb{C})$, on a
\begin{equation}\label{Sylvester}
\det(I_n + A B) = \det(I_d + B A),
\end{equation}
ce qui permet de réduire la dimension du problème, dans le cas général d'une perturbation de rang $d \geq 1$ fixé\footnote{Pour une preuve de cette identité classique, voir Oraux X-ENS \cite{OrauxXENS}.}.
De ces discussions a résulté l'article de Tao \cite{Tao2013}, dont nous présentons ci-dessous quelques résultats emblématiques, pour les perturbations de rang $d=1$.

On considère une matrice $G$ dont les entrées sont des variables aléatoires complexes i.i.d. On peut par exemple prendre des variables gaussiennes complexes: c'est ce qui s'appelle une matrice de Ginibre complexe, en hommage aux premiers résultats emblématiques de Jean Ginibre \cite{Ginibre1964} sur ce modèle. Ainsi donc, considérons, pour fixer les idées:
\begin{equation}\label{Ginibre}
(G_{ij})\ i.i.d., \quad G_{ij} \disteq \mathcal{N}_{\mathbb{C}}(0,\tfrac{1}{N}).
\end{equation}
Les valeurs propres de cette matrice $G$ convergent en distribution vers la loi circulaire de Girko, ce qui signifie que la position d'une valeur propre typique se rapproche en distribution d'un point uniforme sur le disque unité, quand $N\rightarrow \infty$. On sait (voir notamment les travaux de Djalil Chafa{\"{i}} \cite{Chafai2010}) qu'une perturbation de faible rang ne modifie pas la loi circulaire, du moins macroscopiquement, c'est-à-dire que le comportement typique d'une valeur propre n'est pas affecté. Néanmoins, une telle perturbation peut créer des outliers, ce qu'un résultat de Tao confirme, dès lors que la perturbation a elle-même une valeur propre non-nulle assez grande. Il s'agit plus précisément du cas où la perturbation est hermitienne et positive de rang un:
\begin{theorem}[Tao \cite{Tao2013}] Soit $v$ un vecteur unitaire indépendant de $G$ et $\mu$ un  scalaire fixé. Alors le spectre de 
$$ G+\mu\sqrt{N}vv^{*}$$
comporte un unique outlier 
$\lambda_{1} = \mu\sqrt{N}+o(1),$
tandis que le reste du spectre est contenu dans un disque $D(0,1+o(1))$ avec forte probabilité, et converge faiblement vers la loi circulaire pour $N \rightarrow + \infty$. 
\end{theorem}

Remplaçons la valeur $\mu \sqrt{N}$ par un paramètre de temps $t \in \mathbb{R}$, et considérons les trajectoires des valeurs propres des matrices
\begin{displaymath}
G(t) = G + t v v^*, \quad t \in \mathbb{R},
\end{displaymath}
à $N$ fixé. Ces trajectoires sont représentées ci-dessous, où les couleurs ont été réglées de sorte qu'elles tendent vers le bleu clair en $t \rightarrow - \infty$ et vers le bleu foncé en $t \rightarrow +\infty$ en passant par le gris. On constate, tout d'abord, que le spectre limite est le même des deux côtés; c'est un peu surprenant à première vue, mais en fait c'est logique. Pour tout $N$ fixé, les limites des valeurs propres sont les mêmes pour $t\rightarrow\pm\infty$, elles convergent vers des valeurs que l'on peut calculer\footnote{Il s'agit des zéros de la fonction $\mathscr{W}_N$ introduite dans le calcul qui suit.}. Visuellement, les trajectoires des valeurs propres se recollent donc les unes aux autres et forment des boucles dans $\mathbb{C}\cup\{ \infty \}$, la boucle la plus grande étant celle de l'outlier, qui passe par le point à l'infini. \bigskip

\begin{figure}
\begin{center}
\includegraphics[width=0.5\textwidth]{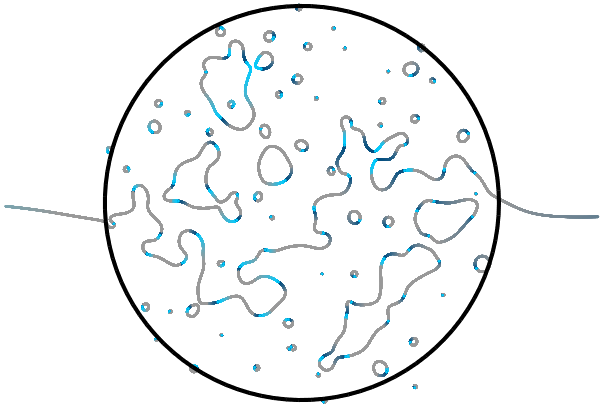}
\end{center}
\caption{Trajectoires de $G+tvv^{*}$ pour une matrice de Ginibre $G$ de taille $100\times 100$.}\label{fig:Gin1}
\end{figure}

Ce résultat de Tao ne fait que prouver rigoureusement une intuition très naturelle; il n'est, à vrai dire, guère surprenant. Mais l'article et les simulations de Rajan et Abbott repèrent également un cas de figure assez inattendu: certaines perturbations de rang un peuvent créer non pas un, mais plusieurs outliers, non fortement séparés au sens défini plus haut. Un autre théorème de Tao vient confirmer et préciser ce phénomène. En voici un énoncé très proche\footnote{La seule différence est que l'article \cite{Tao2013} considère une autre distribution sur le vecteur $w$, plus adaptée au contexte des matrices synaptiques de \cite{RajanAbbott2006}; mais l'esprit est exactement le même.}:
\begin{theorem}\label{Tao2} Soit $v$ un vecteur unitaire fixé, $w$ un vecteur unitaire aléatoire uniforme indépendant de $G$, et $\mu$ un scalaire fixé, alors le spectre de
$$G+\mu\sqrt{N}v w^{*}$$
a un certain nombre (aléatoire) d'ordre $O(1)$ d'outliers, dont la distribution correspond asymptotiquement à celle des zéros d'une série entière gaussienne.
\end{theorem}
Une image des trajectoires dans ce cas de figure est donnée ci-dessous, et la différence avec le cas précédent est assez flagrante. Alors: d'où vient ce nuage de points au-delà du disque unité et pourquoi n'apparaît-il pas lorsque $w=v$? L'heuristique est la suivante\footnote{Cette approche est assez naturelle et s'applique également aux cas hermitiens tels que ceux traités par exemple en \cite{BenaychGuionnetMaida,BenaychRao}.}. On considère le processus aléatoire
$$G(t):= G + t vw^{*},$$
pour des vecteurs $v$ et $w$ que l'on pourra spécifier plus tard.

\begin{figure}[H]
\begin{center}
\includegraphics[width=0.5\textwidth]{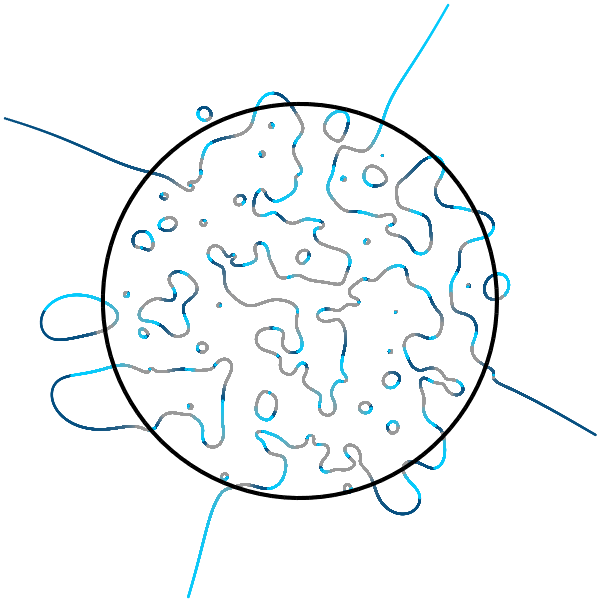} 
\end{center}
\caption{Trajectoires de $G+tvw^{*}$ pour une matrice de Ginibre $G$ de taille $100\times 100$ dans le cas où $w$ est uniforme.}\label{fig:Gin2}
\end{figure}
Pour caractériser le spectre de $G(t)$, on commence par écrire, pour tout $z \in \mathbb{C} \backslash \mathrm{Sp}(G)$ et $t \neq 0$:
\begin{align*}
z\in\text{Sp}(G(t)) & \Leftrightarrow\det(G(t)-z I_N)=0\\
 & \Leftrightarrow\det(I_N+(G-zI_N)^{-1}tvw^{*})=0,
\end{align*}
puis on utilise une propriété élémentaire que les anglo-saxons -- avec la flexibilité syntaxique désarmante dont ils ont le secret -- appellent tout simplement le \textit{matrix-determinant lemma}:
\begin{equation}\label{matrix_det}
\det(I_N+(G-zI_N)^{-1}tvw^{*})
= 1+tw^{*}(G-zI_N)^{-1}v.    
\end{equation}
Remarquons que c'est un cas particulier de l'identité de Sylvester~\eqref{Sylvester}, que l'on peut aussi vérifier ici directement en considérant les valeurs propres de ces opérateurs, puisque $vw^*$ est de rang un. On trouve donc l'équivalence:
\begin{align*}
z\in\text{Sp}(G(t))
 & \Leftrightarrow1+tw^{*}(G-zI_N)^{-1}v=0\\
 & \Leftrightarrow w^{*}(z I_N -G)^{-1}v=\frac{1}{t}.
\end{align*}

Autrement dit le spectre de la matrice perturbée $G(t)$ correspond à des lignes de niveau de la fonction méromorphe aléatoire
\begin{equation}\label{weighted_resolvent}
\mathscr{W}_N(z) := w^{*}(zI_N-G)^{-1}v.
\end{equation}
Cette fonction dépend des variables aléatoires $G$, $v$ et $w$; de manière typique, elle a $N$ pôles qui sont les valeurs propres de $G$. Pour tout $z$ au-delà du rayon spectral (remarquons qu'il suffit de prendre $|z|>1+\epsilon$ pour $\epsilon>0$ fixé\footnote{Conditionner sur l'événement $\Sp (G) \subset D(0,1+\epsilon)$, ne pose pas de problème car on sait que le complémentaire de cet événement a une probabilité exponentiellement faible en $N$: typiquement, une matrice aux entrées i.i.d. n'a aucun outlier, comme cela a été démontré récemment avec brio sous des hypothèses optimales \cite{BordenaveChafaiGZ}.
}), on peut écrire:
\begin{equation}\label{W_N}
\mathscr{W}_N(z) = \frac{w^*v}{z} + \sum_{k \geq 1} \frac{w^* G^k v}{z^{k+1}} .
\end{equation}
On s'intéresse à l'ensemble des outliers au-delà de $1+\epsilon$,
\begin{displaymath}
\left\{ z \in \mathbb{C} \ \ \middle\vert\ \ \ |z|>1+\epsilon, \ \mathscr{W}_N(z) = \frac1t \right\}.
\end{displaymath}
Dans le premier cas de figure ($w=v$), le premier terme de \eqref{W_N} est $z^{-1}$, tandis qu'une analyse assez rapide montre que le reste de l'expression est de taille typique $N^{-1/2}$, en accord avec la normalisation de la matrice $G$ aux entrées i.i.d. \eqref{Ginibre}. La caractérisation du spectre de $G(t)$ devient donc
\begin{equation}\label{approxW}
\frac1t =
\mathscr{W}_N(z) = \frac1z + O(N^{-\frac12})
\end{equation}
et l'on s'attend (avec raison) à un unique outlier dans un voisinage de $z = t$. Ce raisonnement peut être rendu parfaitement rigoureux à l'aide du théorème de Rouché\footnote{On rappelle ici ce théorème incontournable: si $f$ et $g$ sont des fonctions holomorphes sur un domaine $\Omega$ et vérifiant $|f-g| < |g|$ le long d'un lacet simple $\gamma \subset \Omega$, alors $f$ et $g$ ont le même nombre de zéros sur le compact~$K$ qui constitue l'intérieur du lacet. Ce résultat possède des applications très utiles et bien connues des \textit{dog walkers} new-yorkais: si l'on veut éviter que la laisse s'enroule autour d'un arbre, il suffit de faire en sorte qu'elle soit toujours plus courte que la distance qui sépare le \textit{dog walker} de l'arbre. Les mauvaises langues, néanmoins, prétendent que ce principe était déjà connu et appliqué avant l'émergence de l'analyse complexe.} \cite{Rouche}, et en donnant un sens probabiliste précis à l'approximation \eqref{approxW}. 

\bigskip
Dans le second cas de figure (Théorème~\ref{Tao2}) où $w$ est choisi au hasard, uniformément, et indépendamment de $v$, alors $w^*v$ est aussi un terme d'ordre $N^{-1/2}$, et l'on peut prouver par des techniques classiques que les coefficients de \eqref{W_N}, correctement normalisés, convergent vers des gaussiennes indépendantes:
$$
\forall k \geq 0, \ \sqrt{N} w^* G^k v \distconv g_k.
$$
Notons $g$ la série entière gaussienne correspondante:
\begin{equation}\label{g_series}
    g(z)=\sum_{k\geq0} g_{k}z^{k}.
\end{equation}
Ainsi donc, on obtient\footnote{Sans surprise, mais à quelques détails techniques près, résolus par Krishnapur \cite{Krishnapur} et dont nous faisons grâce aux lecteurs.} la convergence
$$
\sqrt{N} \mathscr{W}_N(z) \distconv g\left(\frac1z\right).
$$
Ce qui correspond bien à l'affirmation du Théorème~\ref{Tao2}: avec $t = \mu \sqrt{N}$ pour $\mu$ fixé, les valeurs propres de $G(t)$ correspondent, à la limite, aux inverses des zéros de la fonction $$g(z) - \frac{1}{\mu},$$ qui sont typiquement en nombre fini, d'ordre $1$.

\bigskip
Outre les applications déjà citées (\cite{RajanAbbott2006,Schuessler}) à la modélisation des réseaux de neurones biologiques, des perturbations de faible rang de certaines matrices apparaissent également dans l'étude des réseaux de neurones artificiels et tout particulièrement dans les LLM (Large Language Models) tels que le célèbre GPT-4. On peut mentionner, par exemple, la technique LoRA (Low Rank Approximation \cite{LoRA}) qui permet d'adapter un modèle entraîné au départ pour un certain type de tâche à réaliser une tâche différente à moindre frais; on peut ainsi personnaliser GPT-4 en lui apportant des éléments absents de sa base de donnée initiale. Plus généralement, ces phénomènes d'apparition d'un ou plusieurs outliers au-delà d'une certaine force de la perturbation peuvent être compris comme analogues à l'acquisition d'une certaine faculté au-delà d'une certaine dose d'entraînement; de tels effets de seuils sont observés aussi bien pour un système apprenant à jouer au casse-brique via un algorithme de \textit{reinforcement learning} que pour l'apprentissage des mathématiques par un être humain\footnote{Conformément à l'aphorisme bien connu de Von Neumann: \textit{en mathématiques, on ne comprend pas les choses, on s’y habitue.}}.

\bigskip
Sur le plan purement mathématique et tout particulièrement d'un point de vue probabiliste, la forme de ces trajectoires, la taille et la longueur de leurs concaténations, et tout particulièrement celles de l'unique chemin entre les outliers dans le cas $v=w$, sont autant de paramètres aléatoires que l'on peut légitimement se donner pour défi de quantifier.

\bigskip
Ajoutons en guise d'ouverture que ces trajectoires suivent une équation assez remarquable, que voici:
\begin{equation}\label{EDO}
\lambda_j''(t) = 2 \lambda_j'(t) \sum_{k \neq j} \frac{\lambda_k'(t)}{\lambda_j - \lambda_k}, \quad t>0,
\end{equation}
où les valeurs initiales $\lambda_i(0), \lambda_i'(0)$ sont déterminées par la matrice $G$ et par la perturbation -- autrement dit, tout l'aléa est contenu dans ces conditions initiales, à partir desquelles se produit une évolution déterministe (mais potentiellement chaotique en raison des singularités, lorsque deux valeurs propres se retrouvent proches l'une de l'autre). Pour étudier en détail ces trajectoires, "il suffirait" de réussir à combiner une connaissance assez précise de l'aléa initial avec une étude quantitative de l'équation \eqref{EDO} -- qui après tout n'est jamais qu'une équation différentielle ordinaire... Avis aux amateurs!

\section{Les chaises musicales: perturbation anti-hermitienne d'une matrice aléatoire hermitienne}

\begin{figure}[H]
\begin{center}
\includegraphics[scale=.2]{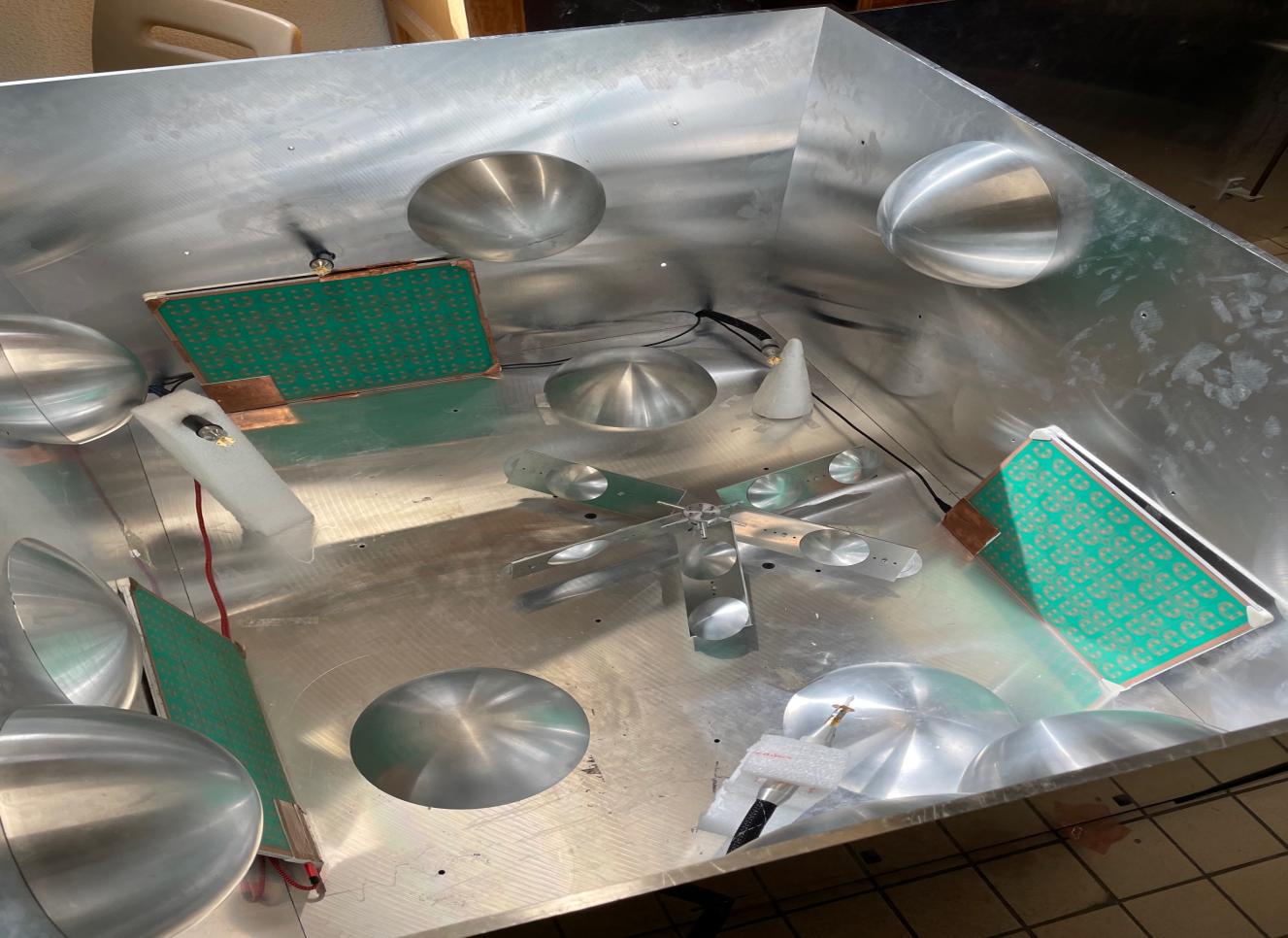}
\end{center}
\caption{Exemple d'une installation expérimentale pour observer des phénomènes chaotiques de dispersion quantique. Source: U.\ Kuhl.}\label{fig:Rev_chamber}
\end{figure}

Les justifications historiques du modèle suivant nous emmènent vers la théorie de la dispersion quantique chaotique (\textit{quantum chaotic scattering}), qui étudie par exemple des expériences faites avec des chambres de réverbérations chaotiques comme celles utilisées à l'Institut de Physique de Nice\footnote{INPHYNI, Université Côte d'Azur} (voir par exemple l'article fondateur \cite{Verbaarschot1985} et les travaux plus récents \cite{GKLM2016,Poli}). Il s'agit d'un appareil encombrant au mode d'emploi compliqué, mais dont le fonctionnement est essentiellement similaire à celui d'un four à micro-ondes: une source émet des ondes électromagnétiques qui sont réfléchies par les murs de la cavité. Mais au lieu de créer une onde stationnaire (qui permet de réchauffer un plat de manière contrôlée), une chambre réverbérante chaotique, comme celle de la figure \ref{fig:Rev_chamber}, doit avoir plusieurs irrégularités qui \textit{dispersent} les ondes d'une manière plus complexe et imprévisible. 

Schématiquement, la dispersion d'une onde électromagnétique est décrite par un opérateur $S$ qui transforme l'onde entrante $\psi_{in}$ en une onde sortante $\psi_{out}$ par l'équation
\begin{equation}\label{eq-inout}
\psi_{out}=S\psi_{in}.
\end{equation}
L'opérateur de dispersion $S$ (\emph{scattering operator}) encode tout le processus et est lié au Hamiltonien $H$, un opérateur auto-adjoint avec un spectre qui correspond aux niveaux d'énergie du système quantique fermé. Cependant, suivre chaque onde de manière déterministe est impossible (par définition, pour un système chaotique) et on utilise plutôt un modèle statistique: c'est là qu'entrent en jeu les matrices aléatoires! Un bon candidat pour modéliser $H$ est une matrice hermitienne aux entrées gaussiennes et essentiellement i.i.d.\ à la symétrie hermitienne près (c'est ce que l'on appelle le GUE, pour \textit{Gaussian Unitary Ensemble}). L'hypothèse de travail est que l'on retrouve le spectre du Hamiltonien, dans la limite où $N \rightarrow \infty$. En réalité, la dispersion n'est pas un processus parfaitement isolé et le système quantique peut échanger de l'énergie avec son environnement. Au lieu d'utiliser la matrice $H$ directement, le système quantique ouvert est donc associé à un Hamiltonien \emph{effectif} $H_{eff}=H+\ri\Gamma$ qui n'est plus auto-adjoint (voir \cite{FyodorovSavin2011} pour les détails). Dans le cas le plus simple, l'interaction avec l'environnement passe par $M$ canaux désignés à l'avance, et l'on a $\Gamma=VV^*$ avec une matrice $V\in\mathbb{C}^{N\times M}$. Ces considérations sont à l'origine de l'intérêt des physiciens (et, à leur suite, des mathématiciens) pour les perturbations anti-hermitiennes de faible rang de matrices hermitiennes. Comme précédemment, le cas le plus simple et le plus connu est celui d'une perturbation anti-hermitienne de rang $1$, pour laquelle un unique outlier apparaît. En physique, l'étude de l'émergence de cet outlier (aussi appelé \textit{broad resonance}) s'inscrit dans la pratique plus générale d'une activité connue sous le nom de capture des résonances (\textit{resonance trapping}).

\bigskip
Soient donc: $H$ une matrice complexe hermitienne, $v$ un vecteur unitaire, $t>0$ un paramètre de temps, et considérons la matrice perturbée
\begin{equation}
    G(t)=H + \ri tvv^{*}.
\end{equation}
On s'intéresse au spectre de $G(t)$, dont les trajectoires sont représentées ci-dessous.

\begin{figure}[H]
\begin{center}
\includegraphics[width=0.8\textwidth]{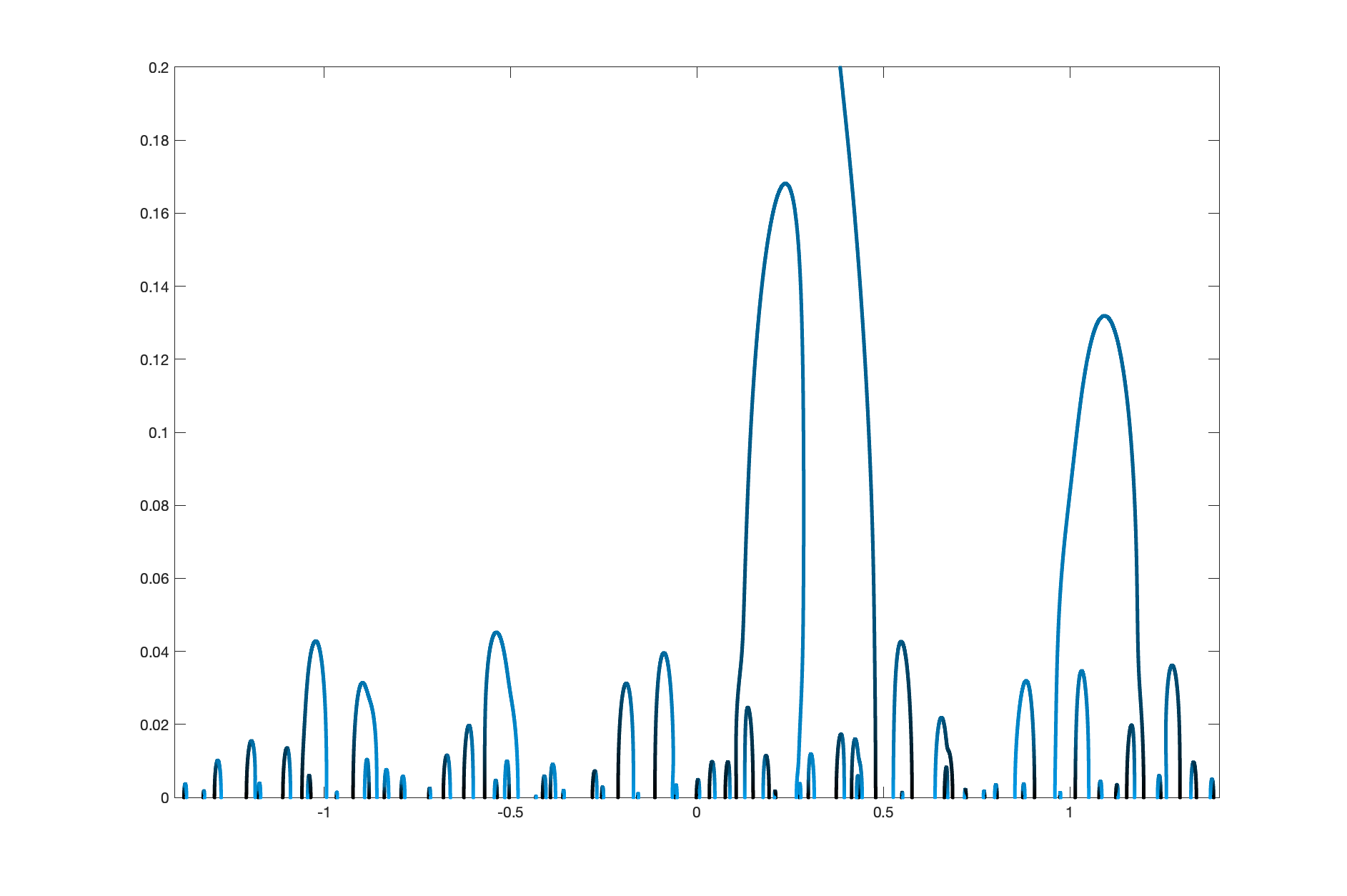} 
\end{center}
\caption{Trajectoires des valeurs propres pour $H$ une matrice GUE de taille $100\times 100$. L'évolution du paramètre~$t$ est indiquée par une variation de couleur{,} du noir ($t{=}0$) au bleu ($t\rightarrow\infty$).}\label{fig:musical_chairs}
\end{figure}

Du côté des physiciens, il est connu depuis longtemps qu'une transition se produit autour de~${t=1}$: voir par exemple les travaux pionniers de Dittes, Heyes et Rotter~\cite{DHR1991} qui font allusion à cette échelle de temps critique; ceci a par la suite été vérifié expérimentalement (voir notamment~\cite{PRSB2000}).

\bigskip
Du côté mathématique, on peut établir les résultats suivants lorsque $H$ est une matrice du GUE: 
\begin{itemize}
\item pour $t=0$, les valeurs propres de $G(0)=H$, sont proches en distribution de la loi du demi-cercle, de densité $\frac{1}{2\pi} \sqrt{4-x^2}$ sur $[-2,2]$. 
\item Pour $t\leq1$, il n'y a pas d'outlier, au sens où l'on peut prouver que toutes les valeurs propres restent proches (distance d'ordre $\frac1N$) de la droite réelle avec forte probabilité.
\item Pour tout $t>1$ fixé, il y a un unique outlier, et fortement séparé; citons par exemple le résultat suivant:
\end{itemize}
\begin{theorem}[O'Rourke et Matchett Wood \cite{ORourkeWood2017}]\label{thm:Wood} Pour $t>1$ et $\epsilon >0$ fixés, avec forte probabilité il y a un outlier $\lambda_{1}$ dans un voisinage de $\ri(t-\frac{1}{t})$. Pour le reste du spectre, on a 
$$ \forall \epsilon >0, \ \forall j \geq 2, \ \Im\lambda_{j}<\frac{N^{\epsilon}}{N},$$ et les parties réelles $(\Re(\la_j))_{j=2}^N$ convergent toujours vers la loi du demi-cercle dans la limite~${N \rightarrow \infty}$.
\end{theorem}
\begin{itemize}
\item enfin, lorsque $t \rightarrow \infty$, l'outlier tend vers l'infini, tandis que les autres valeurs propres convergent vers des valeurs réelles\footnote{Ces limites correspondent aux $N-1$ valeurs propres d'un mineur de $H$. En particulier, elles sont entrelacées avec les $N$ valeurs propres initiales.}.
\end{itemize}
Le comportement général de ces valeurs propres justifie la métaphore des \textit{chaises musicales}, dans la mesure où elles suivent des trajectoires relativement imprévisibles dans le demi-plan supérieur avant de revenir sur la droite réelle, à l'exception d'une valeur propre (l'outlier) qui est exclue du groupe.

\bigskip
Remarquons que pour ce système, la partie imaginaire de la trace est connue: on a 
$$\sum_{j\geq1}\Im(\lambda_{j})=\Im\text{Tr}(H + \ri t vv^*) = t.
$$
Comme par ailleurs, ces parties imaginaires sont positives, l'outlier que l'on attendrait au voisinage de $\ri t$ est en réalité légèrement en-dessous, près d'une position typique $\ri (t-\frac1t)$.
Les parties imaginaires des autres valeurs propres se répartissent les miettes que l'outlier veut bien leur laisser, d'ordre $\frac1t$.

\bigskip
Pour étudier ce système, prenons $z$ dans le demi-plan supérieur ($\Im(z)>0$) et appliquons le \textit{matrix-determinant lemma} comme précédemment. On trouve 
\begin{align*}
z\in\text{Sp}(G(t)) & \Leftrightarrow\det(G(t)-zI_N)=0\\
 & \Leftrightarrow\det(I_N+(H-zI_N)^{-1} \ri t vv^{*})=0\\
 & \Leftrightarrow1+ \ri t v^{*}(H-zI_N)^{-1}v=0\\
 & \Leftrightarrow \mathscr{W}_N(z) := v^{*}(H-zI_N)^{-1}v =\frac{\ri}{t}.
\end{align*}
Jusque là, c'est le même calcul qu'avant, mais à partir de maintenant la méthode doit être adaptée car il serait trop cavalier de développer $v^{*}(H-zI_N)^{-1}v$ dans le régime qui nous intéresse. Heureusement, nous sommes retombés sans y songer sur un objet bien connu. Dans l'immense littérature qui concerne les matrices aléatoires hermitiennes figurent en bonne place les lois locales: le but du jeu est d'estimer la résolvante $(H-zI_N)^{-1}$ lorsque~$z$ est au-dessus du spectre, avec la partie imaginaire la plus petite possible. L'idée générale est que, dans un régime approprié, la trace de la résolvante est proche de la transformée de Stieltjes du demi-cercle:
$$
\frac1N \tr (H-zI_N)^{-1} \simeq m_{\textrm{sc}} (z) := \frac{1}{2\pi}\int_{x=-2}^2  \frac{\sqrt{4-x^2}}{x-z} \dd x.
$$
Alors, bon. Notre fonction $\mathscr{W}_N$ n'est pas une trace de résolvante, certes, mais elle est de la même famille. C'est une sorte de coefficient diagonal de résolvante, qui fait l'objet d'un théorème spécial appelé \textit{loi locale isotrope} \cite{Knowles} et qui permet de justifier l'approximation
$$\mathscr{W}_N(z) \simeq m_{\mathrm{sc}}(z)$$
dans un certain domaine (et avec un terme d'erreur plus grand que pour une loi locale classique: cet objet est plus instable qu'une trace de résolvante). Lorsque ce terme d'erreur est suffisamment petit, néanmoins, on s'attend à ce que l'outlier soit proche d'une solution de l'équation
$$m_{\mathrm{sc}}(z)=\frac{\ri}{t}.$$
Pour une matrice du GUE ou tout autre matrice convergeant vers la loi du demi-cercle, cette prédiction correspond précisément au $\ri (t-\frac1t)$ du Théorème \ref{thm:Wood}. En revanche, notons que selon la même logique, on peut aussi concocter des situations où plusieurs outliers apparaissent pour certaines valeurs de $t$, typiquement si la loi limite n'est pas la loi du demi-cercle, et que la transformée de Stieltjes n'est pas injective \cite{Rochet2017}.

Comme pour le labyrinthe de Ginibre (partie 1), la fonction $\mathscr{W}_N$ caratérise le spectre de $G(t)$ via une relation simple et donc contient en théorie toutes les informations sur les trajectoires~$\lambda_1(t),\dots,\lambda_N(t)$ des valeurs propres. La première chose que l'on peut déduire immédiatement, c'est que les trajectoires ne se croisent (presque) jamais. En effet, en admettant que $\lambda_j(t)=\lambda_k(s)$, on trouve
\begin{equation}
\frac{\ri}{t}=\W_N(\lambda_j(t))=\W_N(\lambda_k(s))=\frac{\ri}{s}
\end{equation}
d'où l'on déduit que $t=s$: c'est un fait déterministe que deux trajectoires ne peuvent se croiser qu'au même temps $t$ ! Mais avec $t=s$, une intersection correspond à un zéro d'ordre au moins deux de $\mathscr{W}_N - \frac{\ri}{t}$ qu'on peut traduire en une condition géometrique\footnote{On peut le résumer par une image: les zéros de $\mathscr{W}_N'$ sont des points aléatoires du plan en nombre fini, et le fait qu'un de ces points tombe pile sur une trajectoire est un événement de probabilité nulle.}, pour conclure enfin par des considérations de dimensions que cet événement a une probabilité nulle.

\bigskip
On peut aussi être plus précis quant à l'échelle de temps au-delà de laquelle le perdant de ce jeu de chaises musicales peut être distingué avec quasi-certitude.

\begin{theorem}[Dubach et Erd\H{o}s \cite{DubachErdos}] Sous des conditions générales\footnote{Ceci couvre le cas du GUE (entrées gaussiennes), mais plus généralement les matrices de Wigner (entrées i.i.d. plus générales, avec des conditions sur les moments).} sur $H$ qui garantissent la loi locale isotrope, pour tout $\epsilon >0$ fixé et $t>1+N^{-1/3+\epsilon}$,
avec forte probabilité, il existe un unique outlier fortement séparé.
\end{theorem}

Pour reprendre l'heuristique en la détaillant un peu : rappelons que l'on a 
\begin{equation}
\mathscr{W}_N(z) := v^{*}(H-zI_N)^{-1}v    
\end{equation}
et que les valeurs propres au temps $t$ sont exactement données par
$$\left\{z \in \mathbb{C}  \ \ \middle\vert\ \ \  \Im z >0, \ \mathscr{W}_N(z)=\frac{\ri}{t} \right\}.$$
La loi locale isotrope de \cite{Knowles} nous donne la domination stochastique\footnote{En matrices aléatoires, on dit que la suite $a_n$ est dominée stochastiquement par la suite $b_n > 0$ (notation: $a_n \prec b_n$) si pour tout $\epsilon>0$, on a $|a_n| < b_n n^{\epsilon} $ avec forte probabilité. C'est un concept pratique pour travailler avec des ordres de grandeur en $n$ pour des variables aléatoires dans la limite $n \rightarrow \infty$.}
\begin{displaymath}
\left| \mathscr{W}_N(z)-m_{\mathrm{sc}}(z) \right| \prec \frac{1}{\sqrt{N \ |\Im z|}}.
\end{displaymath}
Pour tout $t>0, \epsilon>0$, on définit le domaine
\begin{displaymath}
\mathscr{R}_{t,\epsilon}=\left\{\Re(z)^2+\left( \Im(z)-\left(t-\frac1t\right)\right)^{2}>\frac{N^{\epsilon}}{N \cdot \Im(z)} \right\}.
\end{displaymath}
En combinant la domination stochastique avec le théorème de Rouché, on prouve que les
fonctions $w(z)-\frac{\ri}{t}$ et $m_{sc}(z)-\frac{\ri}{t}$ auront
le même nombre de zéros dans ${\mathscr R}_{t,\epsilon}$. En particulier, l'outlier se trouve dans un disque
\begin{displaymath}
D\left(\ri\left(t-\frac1t\right),\frac{N^{\epsilon}}{\sqrt{N \left(t-\frac1t\right)}}\right),
\end{displaymath}
dont le centre est bien la position typique $\ri\left(t-\frac1t\right)$, tandis que le reste des valeurs propres vérifie
\begin{displaymath}
\Im(\lambda_{j})<\frac{N^{\epsilon}}{\left(N\left(t-\frac1t\right)\right)}.
\end{displaymath}
Un calcul direct permet de conclure dès que $t$ est au-delà de l'échelle $1+O(N^{-1/3})$!
Cette technique permet donc d'établir qu'il y a forte séparation de l'outlier au-delà d'une certaine échelle de temps. Néanmoins, elle ne suffit pas à prouver que cette échelle est optimale. Cette confirmation a été fournie peu de temps après par l'article de Fyodorov, Khoruzhenko et Poplavskyi~\cite{FyodorovGUE} qui estime précisément le nombre moyen de points dans un rectangle bien choisi:
\begin{theorem}[Fyodorov, Khoruzhenko et Poplavskyi \cite{FyodorovGUE}] Si l'on définit
\[
{\mathscr N}_{t}(y):=\mathbb{E}[\#\{\lambda\in\text{Sp}(G(t)):\Im(\lambda)>y\}],
\]
on a, quand $N \rightarrow \infty$,
\[
{\mathscr N}_{t}\left ( {y} \right) \
 \sim \ 
\frac{1}{y} \
e^{-Ny(t+\frac{1}{t})} \
I_{1}(2Ny)
\]
 où $I_{1}$ désigne une fonction de Bessel modifée de première espèce et $y$ peut dépendre de $N$.
\end{theorem}
En particulier, pour $t = 1 \pm O(N^{-1/3})$, cette estimation prouve l'existence d'un nuage de valeurs propres dont les parties imaginaires sont d'ordre $O(N^{-1/3})$.
La technique utilisée est un calcul direct assez virtuose, rendu possible par le fait que, dans le cas où $H$ est une matrice du GUE, le système est complètement intégrable, ce qui avait été montré une vingtaine d'années auparavant par deux des trois mêmes auteurs \cite{FyodorovKhoruzhenko}. Dans le cas du GUE, on dispose d'une densité explicite pour la loi jointe des valeurs propres de $H+ \ri tvv^{*}$:
\[
\frac{1}{Z_{N,t}} \delta_{\sum\Im z_{i}=t} \ \prod_{i<j}|z_{i}-z_{j}|^2 \ e^{-\frac{N}{2}(t^{2}+\sum\Re(z_{k}^{2}))}.
\]
par rapport à la mesure de Lebesgue dans le demi-plan supérieur $\{\Im z \geq 0\}$. \\

Nous arrêterons là ce bref survol du modèle des chaises musicales. Mentionnons pour finir une question ouverte très intrigante: puisqu'il s'agit d'un processus continu et sans intersection de trajectoires, de quelle partie du spectre initial (celui de $G(0)=H$) l'outlier provient-il ? Toutes les simulations aussi bien que les calculs tendent à montrer que l'essentiel de l'action se passe dans un petit voisinage de l'origine. Néanmoins, les méthodes actuelles échouent à démontrer rigoureusement que l'outlier provient d'un tel voisinage avec forte probabilité, pour des raisons techniques sérieuses et qui requièrent sans doute que soit inventée une nouvelle approche.

\section{Le modèle UA: perturbation multiplicative d'une matrice unitaire}

Toujours en lien avec la théorie de la dispersion quantique, notre troisième et dernier modèle provient de l'étude des systèmes en temps discret. Il fut introduit par Fyodorov \cite{Fyodorov2001} dans le contexte suivant: considérons un système quantique simple qui interagit avec son environnement. Pour chaque temps (discret) $n$, l'état du système est représenté par un vecteur $\psi_{sys}(n)$, et par ailleurs on considère des signaux $\psi_{in}(n)$ et $\psi_{out}(n)$ qui entrent et sortent du système pour tout~$n$. L'évolution de ce système `ouvert' est alors décrite par la relation
\begin{align*}
\begin{pmatrix}\psi_{sys}(n+1)\\\psi_{out}(n) \end{pmatrix}=\begin{pmatrix} M_{11}&M_{12}\\M_{21}&M_{22}\end{pmatrix}\begin{pmatrix}\psi_{sys}(n)\\\psi_{in}(n) \end{pmatrix}=:M\begin{pmatrix}\psi_{sys}(n)\\\psi_{in}(n) \end{pmatrix}
\end{align*}
assortie d'une contrainte qui garantit la conservation d'énergie. Les éléments de $M$ sont des opérateurs qui décrivent la transformation entre les signaux d'environnement et le système:
\begin{itemize}
    \item $M_{11}$: évolution interne;
    \item $M_{12}$: signal entrant $\rightarrow$ système; \item $M_{21}$: système $\rightarrow$ signal sortant; \item $M_{22}$: transmission entrant $\rightarrow$ sortant.
\end{itemize}
En appliquant quelques transformations \footnote{Nous renvoyons le lecteur à l'article de Fyodorov~\cite{Fyodorov2001} pour les détails techniques que nous omettons ici.}, l'évolution interne prend la forme 
$$M_{11}=U\sqrt{I_N-\tau^*\tau},$$
où $U$ est unitaire et les valeurs propres de $\tau^*\tau$ (pour $\tau$ une matrice rectangulaire $M\times N$) sont associées aux coefficients de transmission qui décrivent le couplage entre le système et son environnement. Dans le cas le plus simple, l'interaction passe par $d$ canaux bien définis, et la perturbation $\tau^*\tau$ prend une forme diagonale: 
\begin{equation}\label{diag_form}
\tau^*\tau=\mathrm{diag}(t_1,\dots,t_d,0,\dots,0).
\end{equation}
Notons que le cas $\tau=0$ correspond à un système quantique \textit{fermé} (sans interactions avec l'environnement), qui est décrit par une évolution unitaire. Pour modéliser sa dynamique chaotique par une approche statistique appropriée, une distribution s'impose: la mesure de Haar sur le groupe unitaire, c'est-à-dire que la matrice $U$ est échantillonnée d'une manière uniforme sur~$U_N(\mathbb{C})$. Dans le jargon historique des matrices aléatoires, on appelle ça le CUE (pour \textit{Circular Unitary Ensemble}).

\begin{figure}[H]
\begin{center}
\includegraphics[width=0.5\textwidth]{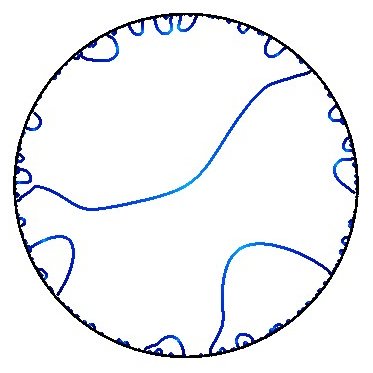} 
\end{center}
\caption{Trajectoires du modèle $UA$ de taille $250\times 250$. L'évolution du paramètre~$t$ est indiquée par la variation de couleur de noir $(|t|{=}1)$ à bleu ($t\rightarrow0$).}\label{fig:UA}
\end{figure}

Voilà pour le contexte. Mathématiquement, on va voir qu'il y a des analogies mais aussi des différences significatives avec les exemples précédents. Pour $U$ une matrice du CUE et
$$A(t) :=(\mathrm{I}_N-(1-t)vv^*)$$
une perturbation de faible rang de l'identité\footnote{Le cas typique est celui d'une matrice diagonale $A(t)=\mathrm{diag}(t,1,\dots,1)$, que l'on obtient avec \eqref{diag_form} quand $d=1$, et que l'on retrouve ici en prenant $v=e_1$. Cette généralisation apanivore permet simplement de retrouver des notations proches de celles des cas précédents.} selon un vecteur unitaire $v$, on considère \emph{le modèle~$UA$} donné par
\begin{equation}\label{eq-UA}
G(t)=UA(t), \qquad t \in [-1,1].  
\end{equation}
On parle d'une \textit{perturbation multiplicative de rang un} pour signifier que la matrice $A$ diffère de l'identité par une matrice de rang un. La dynamique des valeurs propres ressemble à une version circulaire et inclusive\footnote{Il n'y a, en effet, pas de perdant. L'outlier est la valeur propre qui passe par le centre, avant de revenir sur le cercle unité.} des chaises musicales: au départ et à l'arrivée, c'est-à-dire pour $t=\pm1$, $G(t)$ est une matrice unitaire et les valeurs propres sont toutes de module $1$. Au milieu de la dynamique, le spectre de $G(0)$ est identique à celui d'une matrice unitaire tronquée, avec une valeur propre ajoutée à l'origine. Une image de ces trajectoires est donnée sur la figure \ref{fig:UA}.

\bigskip
On observe que la plupart des valeurs propres font de (très) petites
boucles au bord du disque unité, mais il y a une trajectoire qui traverse le disque et passe par zéro (quand $t=0$). C'est cette valeur propre qui est ici l'\textit{outlier} principal: en effet, c'est celle qui s'éloigne le plus du bord du disque, qui est la plus influencée par la perturbation multiplicative, et que l'on peut isoler facilement du reste du spectre. Un argument analogue à celui de la section 2, permet de montrer que les trajectoires de ce modèle $UA$ ne se croisent pas et de déduire un système d'équations différentielles ordinaires pour les décrire (voir~\cite{DubachReker}). Tout cela n'est vrai que pour une perturbation de rang un. Et il y a même une formule pour la densité jointe des valeurs propres de~\eqref{eq-UA} dans la littérature.

\begin{theorem}[Fyodorov \cite{Fyodorov2001}] La densité jointe des valeurs propres du modèle $UA$ est proportionnelle à
\begin{displaymath}
(1-|t|^{2})^{1-N} \ 
\prod_{i<j} |z_{i}-z_{j}|^{2} \
\mathbf{1}_{\{|t|^{2}=\prod|z_{i}|^{2}\}}
\end{displaymath}
par rapport à la mesure de Lebesgue dans le disque unité, avec une constante de normalisation qui est indépendante de $t$.
\end{theorem}

L'échelle de temps $t\sim N^{-1/2}$ joue un rôle privilégié pour le modèle $UA$. En étudiant~\eqref{eq-UA} pour ce choix de $t$, un résultat de Forrester et Ipsen~\cite{ForresterIpsen} montre que la distribution jointe des valeurs propres a comme un air de déjà-vu.

\begin{theorem}[Forrester et Ipsen \cite{ForresterIpsen}] Pour $\mu>0$ fixé et $t=\mu N^{-1/2}$, le spectre de~\eqref{eq-UA} converge en loi vers la distribution des zéros de la série entière gaussienne \eqref{g_series}.
\end{theorem}

Pour s'en convaincre, reprenons le calcul dont nous commençons à avoir pris l'habitude. Pour $|z|<1$ cette fois, on trouve
\begin{align*}
z\in\text{Sp}(G(t)) & \Leftrightarrow\det(G(t)-zI_N)=0\\
 & \Leftrightarrow1-(1-t)v^{*}(U-zI_N)^{-1}Uv=0.
 \end{align*}
Cela nous donne (encore) une fonction méromorphe aléatoire
$$\widetilde{\mathscr{W}_N}(z) := v^{*}(I_N-zU^{*})^{-1}v$$
qui caractérise le spectre de $G(t)$ via la relation
\begin{align*}
z\in\text{Sp}(G(t)) \Leftrightarrow \widetilde{\mathscr{W}_N}(z)=\frac{1}{1-t}.
\end{align*}
Comme $|z|<1$, on peut développer sans crainte:
\begin{align}
\widetilde{\mathscr{W}_N}(z)= \sum_{k\geq0}v^{*}(zU^{*})^{k}v
= 1 + \sum_{k\geq1} c_{k} z^{k}\label{eq-series}
\end{align}
et les coefficients $c_k$, correctement renormalisés par un facteur $\sqrt{N}$ convergent vers des gaussiennes indépendantes. Bref, après ce petit tour de passe-passe, on retombe sur la série gaussienne~\eqref{g_series}. Les valeurs propres au temps $t = \mu N^{-1/2}$ correspondent, à la limite, aux racines de~${g-\mu}$.

\bigskip
Il s'avère que cette échelle $t\sim N^{-1/2}$ considérée par Forrester et Ipsen pour obtenir la convergence ci-dessus régit également l'émergence de l'outlier. C'est en fait à cette échelle que tout se passe.

\begin{theorem}[Dubach et Reker \cite{DubachReker}] Pour tout $\epsilon>0$, $|t|<N^{-1/2-\epsilon}$, il existe un unique outlier fortement séparé. De plus, l'échelle $N^{-1/2}$ est optimale.
\end{theorem}

L'idée est similaire aux arguments exposés dans les sections précédentes: on veut attraper l'outlier au lasso, en utilisant le théorème de Rouché sur un contour approprié. Cependant, la symétrie du modèle $UA$ ne permet pas de travailler avec une approximation complètement déterministe; il faut garder une petite dose d'aléatoire. Pour ce faire, notons que, grâce à~\eqref{eq-series}, la fonction
\begin{displaymath}
\widetilde{\mathscr{W}_N}(z)=v^{*}(I_N-zU^{*})^{-1}v
\end{displaymath}
est donnée par une série et écrivons
\begin{displaymath}
\widetilde{\mathscr{W}_N}(z)=1+v^{*}U^{*}vz+\widetilde{w}_2(z).
\end{displaymath}
Ce sont les fluctuations du terme linéaire qui importent et qu'il convient de garder à part. Comme le CUE est un ensemble assez régulier, l'erreur $\widetilde{w}_2$ est dominée stochastiquement de manière uniforme dans un disque centré en 0 de rayon $1-N^{-\delta}$ pour $\delta>0$. Plus précisément, pour tout $\epsilon>0$, on obtient la borne
\begin{equation}\label{eq-w2bound}
\widetilde{w}_2(z)<\frac{N^\epsilon|z|^2}{\sqrt{N}(1-|z|)}
\end{equation}
avec forte probabilité. Alors, pour
\begin{align*}
f(z)&:=\widetilde{\mathscr{W}_N}-\frac{1}{1-t},\\
g(z)&:=1+v^{*}U^{*}vz-\frac{1}{1-t}
\end{align*} 
on a 
$$|f-g|=|\widetilde{w}_2|$$
et l'on vérifie que $|f-g|<|g|$ sur le domaine
\begin{displaymath}
\widetilde{{\mathscr R}}_{t,\epsilon}=\left \{|z|<1  \ \ \middle\vert\ \ \ \frac{N^{\epsilon}|z|^{2}}{1-|z|}<|v^{*}U^{*}v|\ |z-z_{t}|\right\}
\end{displaymath}
avec forte probabilité, où $z_{t}$ désigne l'unique zéro de $g$. On peut donc appliquer le théorème de Rouché sur le domaine $\widetilde{{\mathscr R}}_{t,\epsilon}$ et prouver que l'outlier est fortement séparé.

\bigskip
Enfin, notons que cet argument s'applique aussi à des matrices aléatoires unitaires au-delà du seul cas du CUE. En fait, c'est vrai pour tout ensemble qui satisfait la borne \eqref{eq-w2bound} pour $\widetilde{\mathscr{W}_N}$ ainsi que l'inégalité 
$$N^{1/2-\epsilon}<v^{*}U^{*}v<N^{1/2+\epsilon}$$
pour tout $\epsilon>0$, avec forte probabilité. Notamment, ce n'est pas une propriété qui dépend de la répulsion des valeurs propres de $U$: elle est toujours vérifiée si ces valeurs propres suivent un processus $\alpha$-déterminantal avec $\alpha\in[-1,1]$. Ces processus fournissent une interpolation continue entre un processus ponctuel avec répulsion (cas \textit{déterminantal}: $\alpha=-1$) et un processus avec attraction (cas \textit{permanental}: $\alpha=1$), en passant par des points indépendants ($\alpha=0$).

Le modèle UA est sans doute celui dans lequel l'espoir de résoudre de nouveaux problèmes par des calculs exacts est le plus concret. Il possède une symétrie complète par rotation qui joue un rôle essentiel dans notre analyse (permettant notamment une preuve directe du caractère optimal de l'échelle $N^{-1/2}$) ainsi que dans le travail à venir de Fyodorov, Khoruzhenko et Prellberg \cite{Fyodorov_CUE}: cette symétrie permet d'invoquer des résultats puissants et parfois contre-intuitifs telles que le théorème de Kostlan (voir \cite{Kostlan, GinibrePowers}) qui affirme que les points d'un processus ponctuel aléatoire en deux dimensions peuvent se repousser fortement tout en ayant des rayons indépendants. Enfin, il semble également naturel de vouloir généraliser les résultats existants à des perturbations de plus haut rang qui préservent la symétrie de rotation. La figure ci-dessous représente les trajectoires pour une perturbation de rang deux.

\bigskip
Ainsi s'achève, avec ces considérations sur les perturbations multiplicatives faiblement non-unitaires, notre bref tour d'horizon des perturbations non-hermitiennes de rang un de matrices aléatoires.

\begin{figure}[H]
\begin{center}
\includegraphics[width=0.5\textwidth]{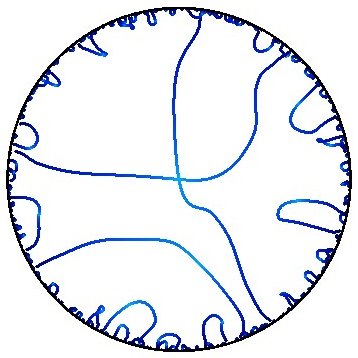}
\end{center}
\caption{En considerant le modèle $UA$ avec une perturbation de rang deux on obtient deux outliers.}\label{fig:UArk2}
\end{figure}

\section*{Épilogue: les vecteurs propres sont-ils de droite ou de gauche?}

Pour une matrice hermitienne $H$ et une valeur propre $
\lambda \in \mathbb{R}$, la tradition veut que les vecteurs propres à droite soient privilégiés, c'est-à-dire que l'on considère par défaut des vecteurs-colonnes $u$ tels que
$$
H u =\lambda u.
$$
Cela dit, si $u$ est un vecteur propre à droite, on peut tout aussi bien constater que le vecteur-ligne $u^*$ est un vecteur propre à gauche, pour la même valeur propre, c'est-à-dire que
$$
u^* H = \lambda u^*.
$$
Autrement dit, pour un vecteur propre dans le monde hermitien, être à gauche ou à droite, ça ne change pas grand chose. Il en va tout autrement dans le monde non-hermitien... Si $G$ est une matrice complexe diagonalisable (ce qui est le cas avec probabilité $1$ pour les modèles habituels de matrices aléatoires), on peut écrire
\begin{displaymath}
G = P \Delta P^{-1}, 
\quad \Delta = \mathrm{diag}(\lambda_1, \dots, \lambda_N),
\end{displaymath}
puis considérer $R_1, \dots, R_N$ les colonnes de $P$ et $L_1, \dots, L_N$ les lignes de $P^{-1}$, qui vérifient
\begin{displaymath}
G R_j = \lambda_j R_j, \quad L_j G = \lambda_j L_j.
\end{displaymath}
On voit donc que $R_i$ est un vecteur propre à droite, et $L_i$ est un vecteur propre à gauche, pour la même valeur\footnote{Il est bon de remarquer que les vecteurs propres, qu'il soient à gauche ou à droite, partagent les mêmes valeurs.} propre $\lambda_i$. Ces deux vecteurs peuvent être très différents... Par ailleurs, aucune de ces deux familles n'est, en général, une famille orthogonale -- mais l'orthonormalité ($u_i^* u_j = \delta_{ij}$) du cas hermitien est remplacée par la bi-orthogonalité. En effet,
\begin{displaymath}
L_i R_j = \delta_{ij},
\end{displaymath}
ce qui découle directement de $P^{-1}P=I_N$.
Une perturbation de la matrice non-hermitienne $G$ modifie tous ces éléments ensemble, de manière couplée. Dans tous les modèles présentés ci-dessus, la dérivée des valeurs propres s'exprime en fonction des valeurs propres \textit{et} des vecteurs propres, et ainsi de suite. Comment résumer cette information pour que ce système titanesque devienne (au moins partiellement) intelligible? Plusieurs méthodes peuvent être envisagées. Une piste privilégiée est l'étude de quantités homogènes dépendant des deux types de vecteurs propres, telles que la matrice $\mathscr{O}$ des overlaps \cite{FyodorovMehlig}
dont les entrées sont données par un produit de produits scalaires hermitiens entre vecteurs propres à droite et à gauche:
$$
\mathscr{O}_{ij} = (L_i L_j^*) (R_j^* R_i).
$$
Mais ceci est une autre histoire...

\bigskip
\textbf{Remerciements:} Cet article a pour point de départ les travaux récents des deux auteurs et un mini-cours donné par G. D. à l'Institut Henri Poincaré dans le cadre du séminaire Matrices Et Graphes Aléatoires (MEGA). Nous tenons à remercier les organisateurs et les participants du~MEGA, tout particulièrement Raphaël Ducatez dont les notes ont été le point de départ de cet article, ainsi que Margaret Bilu, Victor Dubach, Aniss Farès et George-Ioan Stoica pour leurs nombreuses suggestions. Nous avons, enfin, grandement bénéficié de la relecture attentive de deux excellents rapporteurs, Pauline Lafitte et Maxime Février, que nous remercions chaleureusement~!

\renewcommand*{\bibname}{References}

\let\oldthebibliography\thebibliography
\let\endoldthebibliography\endthebibliography
\renewenvironment{thebibliography}[1]{
  \begin{oldthebibliography}{#1}
    \setlength{\itemsep}{0.5em}
    \setlength{\parskip}{0em}
}
{
  \end{oldthebibliography}
}

\bibliographystyle{plain}
\bibliography{refs}

\end{document}